\documentclass{proc-l}
\usepackage{amssymb, graphicx}
\usepackage[cmtip,all]{xy}

\newtheorem{theorem}{Theorem}[section]

\theoremstyle{definition}
\newtheorem{definition}[theorem]{Definition}
\newtheorem{proposition}[theorem]{Proposition}
\newtheorem{example}[theorem]{Example}

\theoremstyle{remark}
\newtheorem{remark}[theorem]{Remark}

\numberwithin{equation}{section}

\begin{document}

\title{Integrability Criterion for Abelian Extensions of Lie Groups}
\author{Pedram Hekmati}
\address{School of Mathematical Sciences, University of Adelaide, 
Adelaide, SA 5005, Australia}
\email{pedram.hekmati@adelaide.edu.au}
\thanks{}


\dedicatory{}


\begin{abstract}
		We establish a criterion for when an abelian extension of infinite-dimensional 
		Lie algebras $\mathfrak{\hat{g}} = \mathfrak{g} \oplus_\omega \mathfrak{a}$ integrates to a 
		corresponding Lie group extension $A \hookrightarrow \widehat{G} \twoheadrightarrow G$, 
		where $G$ is connected, simply connected and $A \cong \mathfrak{a} \slash \Gamma$ for some 
		discrete subgroup $\Gamma \subseteq \mathfrak{a}$. When $\pi_1(G)\neq 0$, the kernel $A$ is 
		replaced by a central extension $\widehat{A}$ of $\pi_1(G)$ by $A$.
\end{abstract}

\maketitle

\section{Introduction}
Given a group $G$ with a normal subgroup $N$, one may construct the quotient group $H = G/N$. 
The theory of group extensions addresses the converse problem. Starting with $H$ and $N$, what different 
groups $G$ can arise containing $N$ as a normal subgroup such that $H \cong G/N$? The problem can be 
formulated for infinite-dimensional Lie groups, but the situation is more delicate. Many familiar 
theorems break down and one must take into account topological obstructions. In particular, Lie's third 
theorem no longer holds and the question of integrability, i.e. whether a Lie algebra corresponds to a 
Lie group, becomes relevant \cite{G}.

The aim of this paper is to establish an integrability criterion for abelian extensions of 
infinite-dimensional Lie groups by generalizing a geometric construction for gauge groups
\cite{LMNS,M2,M3,M4}. A similar idea is employed in \cite{PS} to construct a prequantum bundle for a closed 
2-form with integral periods. For an alternative approach to this problem see \cite{N}. In Sections 
2 and 3 we review the basic definitions of infinite-dimensional Lie groups and their abelian extensions. 
Section 4 gives a detailed account of the construction leading up to the integrability criterion.

\section{Infinite-dimensional Lie groups}
We define infinite-dimensional Lie groups along the lines of \cite{M5}, which should be 
consulted for further details and for concrete examples. The first step is to define the concept of an 
infinite-dimensional smooth manifold. Here the bottom line is to replace $\mathbb{R}^n$ 
(or $\mathbb{C}^n$) by a more general model space on which a meaningful differential calculus can be 
developed. Essentially all familiar constructions in finite dimensions then carry over to the 
infinite-dimensional setting. We consider sequentially complete locally convex topological vector 
spaces. These spaces have the property that every continuous path has a Riemann integral. We adopt the 
following notion of smoothness.
\begin{definition}
Let $E, F$ be sequentially complete locally convex topological vector spaces over $\mathbb{R}$ 
(or $\mathbb{C}$) and let $f:U \to F$ be a continuous map on an open subset $U\subseteq E$. Then 
$f$ is said to be differentiable at $x \in U$ if the directional derivative 
\begin{equation*} df(x)(v) = \lim_{t \to 0} \frac{1}{t}(f(x+tv)-f(x))\end{equation*}
exists for all $v \in E$. It is of class $C^{1}$ if it is differentiable at all points of $U$ and
\begin{equation*} df: U\times E \to F,  \; (x,v) \mapsto df(x)(v) \end{equation*}
is a continuous map on $U \times E$. Inductively we say that
$f$ is of class $C^{n}$ if $df$ is a map of class $C^{n-1}$ and of
class $C^{\infty}$ or smooth if it is of class $C^{n}$ for all $ n \geq 1$.
\end{definition}
This definition coincides with the alternative notion of \textit{convenient }smoothness \cite{K2} 
on Fr\'echet manifolds. A smooth manifold modeled on a sequentially complete locally convex topological 
vector space $E$ is a Hausdorff topological space $M$ with an atlas of local charts $\{(U_i, \phi_i)\}$ 
such that the transition functions $\phi_i \circ \phi_j^{-1}:\phi_j(U_i \cap U_j)\to \phi_i(U_i \cap U_j)$ are smooth on overlaps.
A Lie group $G$ is a smooth manifold endowed with a group structure such that the operations of 
multiplication and inversion are smooth. The Lie algebra $\mathfrak{g}$ is defined as the space of 
left-invariant vector fields. A vector field $X:G \to TG$ is left-invariant if 
\begin{equation*} 
	L_{g*}X = X, \; \; \forall g \in G \ 
\end{equation*}
where $L_{g*}$ denotes the pushforward map induced by the diffeomorphism $L_g:G \to G, \; h \mapsto gh$. 
By definition, $X$ is completely determined by its value at the identity and $\mathfrak{g}$ is therefore 
identified with $T_{\mathbf{1}}G$ as topological vector spaces endowed with the continuous Lie bracket 
of vector fields. The most striking feature of infinite-dimensional Lie theory is that results on 
existence and uniqueness of ordinary differential equations and the implicit function theorem cease to
hold in general beyond Banach Lie groups. Therefore a priori there is no exponential map and even if it exists, 
it does not have to be locally bijective. The existence and smoothness of the exponential function 
hinges on the notion of regularity.
\begin{definition}
A Lie group $G$ is called \textit{regular} if for each $X \in C^{\infty}([0,1],\mathfrak{g})$, there exists
$\gamma \in C^{\infty}([0,1],G)$ such that 
\begin{equation*}
	\gamma ' (t) = L_{\gamma(t)*}(\mathbf{1}).X(t), \ \ \ \ \ \gamma(0)=\mathbf{1}
\end{equation*}
and the evolution map
\begin{equation*}
	{\rm evol}_G:C^{\infty}([0,1],\mathfrak{g}) \to G, \ \ \ \ \ X \mapsto \gamma(1)
\end{equation*}
is smooth.
\end{definition}
In other words every smooth curve in the Lie algebra should arise, in a smooth way, as the left 
logarithmic derivative of a smooth curve in the Lie group. Note that regularity is a sharper condition 
than the requirement that the exponential map should be defined and smooth. Indeed if $\gamma(t)$ is 
the curve corresponding to the constant path $X(t) = X_0$ for some $X_0\in \mathfrak{g}$, then 
$\gamma(1) = {\rm exp}(X_0)$. All known Lie groups modeled on sequentially complete locally convex 
topological vector spaces are regular \cite{G}. In the convenient setting for calculus, it has been 
shown \cite{M1} that all connected regular abelian Lie groups are of the form $\mathfrak{a}/\Gamma$ 
for some discrete subgroup $\Gamma \subseteq \mathfrak{a}$ of an abelian Lie algebra $\mathfrak{a}$. 
Moreover, parallel transport exists for connections on principal bundles with regular structure 
group \cite{K1}. Important examples of regular Lie groups include gauge groups $C^\infty(M,G)$ and 
diffeomorphism groups ${\rm Diff}(M)$, where $M$ is a smooth compact manifold and $G$ is a 
finite-dimensional Lie group. A Lie group is called locally exponential if the exponential function 
exists and is a local diffeomorphism at the identity. Gauge groups have this property, but it does 
not hold true for diffeomorphism groups.

We digress to say a few words about the cohomology of Lie groups and Lie algebras.

\subsection{Lie group cohomology} Let $G$ be a Lie group. An abelian Lie group $A$ is called a smooth 
$G$-module if there is a smooth $G$-action on $A$ by automorphisms $G\times A \to A, (g,a) \mapsto g.a$. 
The set of smooth maps $f:G^n \to A$ such that $f(g_1,\dots,g_{n})=0$ whenever 
$g_j = \textbf{1}$ for some $j$ are called \textit{n-cochains} and form an abelian group $C^n(G,A)$ under 
pointwise addition. A cochain complex
\begin{equation*}
	\dots \to C^{n-1}(G,A) \xrightarrow{\delta_{n-1}} C^n(G,A) \xrightarrow{\delta_n} C^{n+1}(G,A) \to \dots
\end{equation*}
is generated by the homomorphisms $\delta_n:C^n(G,A) \to C^{n+1}(G,A)$ defined by
\begin{eqnarray*}
	(\delta_n f)(g_1,\dots,g_{n+1}) &=& g_1 . f(g_2,\dots,g_{n+1}) + 
	\sum_{i=1}^n(-1)^if(g_1,\dots,g_ig_{i+1},\dots,g_{n+1}) \\ && + (-1)^{n+1}f(g_1,\dots,g_n) \ 
\end{eqnarray*}
and satisfying $\delta_{n} \circ \delta_{n-1} =0$. 
Let $Z^n(G,A)={\rm ker } \ \delta_n$ and $B^n(G,A)={\rm im \; }\delta_{n-1}$ denote the subgroups of 
\textit{n-cocycles} and \textit{n-coboundaries} respectively. The $n$-th Lie cohomology group is defined by  
\begin{equation*}
	H^n(G,A) = \frac{Z^n(G,A)}{B^n(G,A)}\ .
\end{equation*}
\\

\subsection{Lie algebra cohomology}
Let $\mathfrak{g}$ and $\mathfrak{a}$ be topological Lie algebras. Then $\mathfrak{a}$ is a continuous 
$\mathfrak{g}$-module if it is abelian and there is a continuous 
$\mathfrak{g}$-action, $\mathfrak{g} \times \mathfrak{a} \to \mathfrak{a}, (X,v) \mapsto X.v$. Denote by 
$C^n(\mathfrak{g},\mathfrak{a})$ the vector space of continuous alternating multilinear maps 
$\omega:\mathfrak{g}^n \to \mathfrak{a}$. A cochain complex
\begin{equation*}
	\dots\to C^{n-1}(\mathfrak{g},\mathfrak{a}) \xrightarrow{d_{n-1}} C^n(\mathfrak{g},\mathfrak{a})
 \xrightarrow{d_n} C^{n+1}(\mathfrak{g},\mathfrak{a}) \to \dots
\end{equation*}
is generated by the linear maps $d_n:C^n(\mathfrak{g},\mathfrak{a}) \to C^{n+1}(\mathfrak{g},\mathfrak{a})$ 
given by Palais' formula
\begin{eqnarray*} 
	(d_n \omega)(X_1,\dots,X_{n+1}) &=& \sum_{i=1}^{n+1} (-1)^{i+1}X_i .
\omega(X_1,\stackrel{i}{\check{\cdots}},X_{n+1}) \\  && + \sum_{i<j}(-1)^{i+j}
\omega([X_i,X_j],X_1,\stackrel{i}{\check{\cdots}}\stackrel{j}{\check{\cdots}},X_{n+1})
\end{eqnarray*}
and satisfying $d_{n}\circ d_{n-1} =0$. 
Let $Z^n(\mathfrak{g},\mathfrak{a})={\rm ker } \ d_n$ and $B^n(\mathfrak{g},\mathfrak{a})={\rm im \; }d_{n-1}$ 
denote the subspaces of \textit{n-cocycles} and \textit{n-coboundaries} respectively. The $n$-th Lie 
cohomology group is given by the quotient space
\begin{equation*}
	H^n(\mathfrak{g},\mathfrak{a}) = \frac{Z^n(\mathfrak{g},\mathfrak{a})}{B^n(\mathfrak{g},\mathfrak{a})} \ .
\end{equation*}

In the context of Lie group and Lie algebra extensions, the second cohomology group is important in classifying 
topologically trivial abelian extensions as we will see. For $n\geq 1$ there is a `derivation map' 
$D_n: H^n(G,A) \to H^n(\mathfrak{g},\mathfrak{a})$ given by \cite{N} 
\begin{equation*} 
(D_n f)(X_1,\dots,X_n) = \frac{\partial^n}
{\partial t_1 \dots \partial t_n} \sum_{\sigma \in S_n} {\rm sgn}(\sigma)
f\left(\gamma_{\sigma(1)}(t_{\sigma(1)}), \dots, \gamma_{\sigma(n)}(t_{\sigma(n)})\right)\Big|_{t_i = 0} 
\end{equation*}
where $\gamma_1(t_1),\dots,\gamma_n(t_n)$ is any set of smooth curves in $G$ satisfying 
$\gamma_i(0) = \mathbf{1}$ and $\gamma_i'(0) = X_i \in \mathfrak{g}$.

\section{Abelian extensions}
\begin{definition}
An extension of Lie groups is a short exact sequence with smooth homomorphisms
\begin{equation*}
	\mathbf{1} \to A \xrightarrow{i} \widehat{G} \xrightarrow{p} G \to \mathbf{1}
\end{equation*}
such that $p$ admits a smooth local section $\sigma: U \to \widehat{G}$, $p \circ \sigma = {\rm id}_U$, 
where $U \subset G$ is an open identity neighborhood.
\end{definition}
The existence of a smooth local section means that $\widehat{G}$ is a principal $A$-bundle over $G$. The extension 
is called \textit{abelian} if $A$ is abelian and \textit{central} if $i(A)$ lies in the center $Z(\widehat{G})$. 
Two extensions $\widehat{G}_1$ and $\widehat{G}_2$ are \textit{equivalent} if there exists a smooth homomorphism
$\phi:\widehat{G}_1\to\widehat{G}_2$ such that the following diagram commutes:
\begin{displaymath} 
\xymatrix{ 
{A}\ar[r]^{i_1}\ar[d]_{id_A} 
&{\widehat{G}_1}\ar[r]^{p_1}\ar[d]^{\phi}&{G}\ar[d]^{id_G}\\ 
{A}\ar[r]^{i_2} &{\widehat{G}_2}\ar[r]^{p_2}&{G}} 
\end{displaymath}
It is straightforward to verify that $\phi$ must be a Lie group isomorphism. The definition
for Lie algebras is analogous. 
\begin{definition}
An extension of topological Lie algebras is a short exact sequence with continuous homomorphisms
\begin{equation*}
	\mathbf{0} \to \mathfrak{a} \xrightarrow{i} \mathfrak{\hat{g}} \xrightarrow{p} \mathfrak{g} \to \mathbf{0} \ .
\end{equation*}
\end{definition}
Two extensions $\mathfrak{\hat{g}}_1$ and $\mathfrak{\hat{g}}_2$ are said to be \textit{equivalent} if there is an 
isomorphism of topological Lie algebras $\phi: \mathfrak{\hat{g}}_1 \to \mathfrak{\hat{g}}_2$ such that the following 
diagram commutes:
\begin{displaymath} 
\xymatrix{ 
{\mathfrak{a}}\ar[r]^{i_1}\ar[d]_{id_\mathfrak{a}} 
&{\widehat{\mathfrak{g}}_1}\ar[r]^{p_1}\ar[d]^{\phi}&{\mathfrak{g}}\ar[d]^{id_\mathfrak{g}}\\ 
{\mathfrak{a}}\ar[r]^{i_2} &{\widehat{\mathfrak{g}}_2}\ar[r]^{p_2}&{\mathfrak{g}}} 
\end{displaymath}
Next we show how abelian extensions can be constructed explicitly. We will assume that as a 
principal bundle, $\widehat{G}$ is smoothly trivial; i.e. there exists a smooth global section 
$\sigma: G\to \widehat{G}$. 
\begin{proposition} 
	Let $G$ be a Lie group, $A$ a smooth $G$-module and $f \in Z^2(G,A)$ a 2-cocycle. The smooth manifold 
	$G \times A$ endowed with the multiplication
	\begin{equation*}
		(g_1,a_1)(g_2,a_2) = (g_1g_2, a_1 +g_1.a_2 + f(g_1,g_2))
	\end{equation*}
	defines an abelian extension $\widehat{G}= G \times_f A$ of $G$ by $A$.
\end{proposition}
Associativity of the group law follows by the 2-cocycle property. The unit element is $(\mathbf{1},0)$ and 
$(g,a)^{-1} = \left(g^{-1},-g^{-1}.(a+f(g,g^{-1}))\right)$. The extension is smoothly trivial by construction, 
and the conjugation action of $\widehat{G}$ on $A$ induces the smooth $G$-action. When the cocycle is a coboundary, 
the extension is isomorphic to the semidirect product $G \ltimes A$. There is a similar cocycle 
construction for Lie algebras. 
\begin{proposition}
	Let $\mathfrak{g}$ be a topological Lie algebra, $\mathfrak{a}$ a continuous $\mathfrak{g}$-module and 
	$\omega \in Z^2(\mathfrak{g},\mathfrak{a})$ a 2-cocycle. The topological vector space 
	$\mathfrak{g}\oplus \mathfrak{a}$ endowed with the continuous Lie bracket
\begin{equation*}
	[(X_1,v_1),(X_2,v_2)] = ([X_1,X_2], X_1.v_2 - X_2.v_1 + \omega(X_1,X_2))
\end{equation*}
defines a topologically split abelian extension $\hat{\mathfrak{g}}= \mathfrak{g}\oplus_\omega \mathfrak{a}$ 
of $\mathfrak{g}$ by $\mathfrak{a}$. 
\end{proposition}

It turns out that all smoothly trivial abelian extensions arise in this way. Furthermore, two such extensions are equivalent if and only if the 2-cocycles differ by a 2-coboundary \cite{N}.\footnote{ In \cite{N}, the author considers a cohomology theory $H^\bullet_s(G,A)$ 
based on locally smooth cocycles. The cohomology groups $H^n(G,A)$ in the present paper embed as subgroups in $H^n_s(G,A)$; cf. Remark 8.5 in \cite{N}.} The second cohomology groups $H^2(G,A)$ and 
$H^2(\mathfrak{g},\mathfrak{a})$ therefore parametrize the set of equivalence classes of these extensions. The Lie 
algebra of $\widehat{G}= G \times_f A$ is as one would expect $\hat{\mathfrak{g}}= \mathfrak{g}\oplus_{D_2 f}\mathfrak{a}$. 

\begin{remark} An abelian extension $\widehat{G}$ is central if and only if $A$ is a trivial $G$-module. Indeed, since the conjugation action of $\widehat{G}$ on $A$ induces the smooth $G$-action, there is for all $g\in G$, $a\in A$ an element $\hat g \in \widehat{G}$ such that $i(g.a)=\hat g i(a) \hat g^{-1}$. If $i(A)\subset Z(\widehat G)$, then $i(g.a)=i(a)$ and the injectivity of the inclusion map $i$ implies that $g.a = a$, so the $G$-action is trivial. Conversely, if $G$ acts trivially on $A$, then it follows by $i(a)=\hat g i(a) \hat g^{-1}$ that $i(A)$ belongs to the center of $\widehat G$. 
	
In particular, the extension $\widehat{G}$ is central whenever $G$ is connected and the automorphism group $\text{Aut}(A)$ is discrete. The latter occurs for instance when $A$ is a finite-dimensional real connected compact abelian Lie group and hence isomorphic to the torus $\mathbb{T}^n$. 
Interesting examples of abelian extensions which are non-central arise when $A$ is infinite-dimensional. In \cite{M2} Mickelsson describes abelian extensions of gauge groups $C^\infty(M,G)$ by the group of circle valued functions $Map(\mathcal{A},\mathbb{T})$, where $\mathcal{A}$ denotes the affine space of $\mathfrak{g}$-valued 1-forms on a compact smooth manifold $M$. 

Another interesting case is provided by Neeb, Example 9.16 in \cite{N}: Let $Z$ be an abelian group with the Lie algebra $\mathfrak z$. Consider a smooth $Z$-principal bundle $P$ over the compact manifold $M$ together with a connection form $\theta \in \Omega^1(P,\mathfrak z)$. The corresponding curvature form $\omega\in \Omega^2(M,\mathfrak z)$ then defines a 2-cocycle on the Lie algebra of vector fields $\mathcal{V}(M)$. Furthermore, vector fields on $M$ act non-trivially on the commutative gauge algebra $C^\infty(M,\mathfrak z)$ by derivations. The corresponding Lie algebra extension $\widehat{\mathfrak{g}}=\mathcal{V}(M)\oplus_\omega C^\infty(M,\mathfrak z)$ is naturally identified with the Lie algebra $\mathcal{V}(P)^Z$ of $Z$-invariant vector fields on $P$ and lifts to a non-trivial abelian Lie group extension $C^\infty(M,Z) \hookrightarrow \text{Aut}_Z(P) \twoheadrightarrow \text{Diff}(M)_0$.
\end{remark}

\section{Integrability criterion}
In this section we elucidate when an abelian extension of Lie algebras 
$\mathfrak{\hat{g}} = \mathfrak{g} \oplus_\omega \mathfrak{a}$ corresponds to a Lie group extension $\widehat{G}$. 
If $\omega = D_2 f$ for some cocycle $f \in Z^2(G,A)$, then by the previous section a corresponding Lie group extension is 
$\widehat{G}= G \times_f A$. In the general case, $\omega$ must satisfy a certain integrality condition that will 
become apparent by the following construction. The basic idea is to construct $\widehat{G}$ as the quotient of a 
larger group $\mathcal{P}G \times_\gamma A$. This means that in general the extension will be topologically twisted 
and therefore the group multiplication cannot be described by a smooth global 2-cocycle.

Let $G$ be a connected Lie group of $\mathfrak{g}$ and let $A$ be a smooth $G$-module of the form $\mathfrak{a}/\Gamma$ 
for some discrete subgroup $\Gamma \subseteq \mathfrak{a}$. We write $e:\mathfrak{a} \to A$ for the exponential 
(quotient) map and employ a multiplicative notation. The Lie algebra cocycle $\omega \in Z^2(\mathfrak{g},\mathfrak{a})$ 
defines a closed $G$-equivariant 2-form $\omega^{eq} \in \Omega^2(G,\mathfrak{a})$ by
\begin{equation*}
	\omega^{eq}(g)(L_{g*}X, L_{g*}Y) = (L_{g}^*\omega^{eq})(\mathbf{1})(X,Y) = g.\omega(X,Y) \ \ \ \ \ \forall X,Y \in 
	\mathfrak{g} \ .
\end{equation*}
For central extensions, the $G$-action on $A$ is trivial and this is simply the associated left-invariant 2-form.  
Let $\mathcal{P}G$ denote the space of smooth based paths $\hat{g}:[0,1]\to G$ originating at the identity 
$\hat{g}(0) = \mathbf{1}$ and with some arbitrary endpoint $\hat{g}(1) = g$ in $G$. Given the 
$C^\infty$-topology of uniform convergence of the paths and all their derivatives, $\mathcal{P}G$ becomes a smoothly contractible 
Fr\'echet Lie group under pointwise multiplication. It is further a locally trivial $\Omega G$-principal bundle over $G$, where 
$\Omega G$ is the group of smooth based loops, i.e. smooth paths whose endpoints coincide at the identity \cite{CM}. 
Consider $\mathcal{P}G\times A$ and introduce an equivalence relation
\begin{equation*}
	(\hat{g}_1, a_1) \sim (\hat{g}_2, a_1e^{\int_{\pi[\hat{g}_1,\hat{g}_2]} \omega^{eq}})
\end{equation*}
whenever two paths with the same endpoint $\hat{g}_1(1) = \hat{g}_2(1)$ are homotopic and hence form a null homotopic loop. 
There is then a well-defined 2-dimensional surface $\pi = \pi[\hat{g}_1,\hat{g}_2]$ in $G$ bounded by these paths. 
However, the surface $\pi$ is not unique; it depends on the choice of a smooth homotopy $F:[0,1]\times[0,1] \to G$ with 
$F(0,t)= \mathbf{1}$ and $F(1,t)=(\hat{g}_1 * \hat{g}_2)(t)$, where 
\begin{equation*}
	(\hat{g}_1 * \hat{g}_2)(t) = \begin{cases} \hat{g}_1(2t) & \ 0\leq t\leq \frac 1 2 \\ 
	\hat{g}_2(2-2t) & \frac 1 2 \leq t\leq 1 \ .\end{cases}
\end{equation*}
By triangulation $\pi$ can be described by a smooth singular 2-chain. If $\pi'$ is another smooth 2-chain with the 
same boundary, then
\begin{equation*}
	e^{\int_{\pi'} \omega^{eq}} = e^{\int_{\pi'} \omega^{eq}+ \int_{\pi+ \pi^{-}} \omega^{eq}} = 
	e^{\int_{\pi' + \pi^{-}} \omega^{eq}}e^{\int_\pi \omega^{eq}}
\end{equation*}
where $\pi^{-}$ denotes $\pi$ with the opposite orientation. Here $\pi' + \pi^{-}$ is the smooth spherical 2-cycle 
corresponding to the closed surface obtained by gluing together $\pi'$ and $\pi^{-}$ along their common boundary. 
For the equivalence relation to be well-defined we require that $e^{\int_{\pi'+ \pi^-}\omega^{eq}}= \mathbf{1}_A$
or equivalently that 
\begin{equation}\label{condition}
	\int_{c} \omega^{eq} \in \Gamma
\end{equation}
for all spherical cycles $c \in Z_2(G)$ \cite{E}. If $c = \partial b \in B_2(G)$ is a 2-boundary, then this is automatically 
satisfied by Stoke's theorem
\begin{equation*}
	\int_{\partial b} \omega^{eq} = \int_b d\omega^{eq} = 0 
\end{equation*}
and therefore the condition factors through to homology.

We can now proceed to construct the Lie group extension $\hat{G}$ corresponding to $\mathfrak{\hat{g}} = 
\mathfrak{g} \oplus_\omega \mathfrak{a}$ by defining a multiplication on $\mathcal{P}G\times A / \sim $ 
\begin{equation*}
	[(\hat{g}_1,a_1)][(\hat{g}_2,a_2)] = [(\hat{g}_1\hat{g}_2,a_1(\hat{g}_1.a_2)e^{\gamma(\hat{g}_1,\hat{g}_2)})]
\end{equation*}
where $\hat{g}.a := \hat{g}(1).a = g.a$ is the given $G$-action on $A$ and 
$\gamma : \mathcal{P}G \times \mathcal{P}G \to \mathfrak{a}$ is a smooth 2-cocycle. The latter must be defined in 
such a way that it yields the correct Lie algebra cocycle $\omega$ and is compatible with the equivalence relation. 
This is accomplished by choosing
\begin{equation*}
	\gamma(\hat{g}_1,\hat{g}_2) = \int_\sigma \omega^{eq}
\end{equation*}
where $\sigma: \Delta_2 \to G, (t,s) \mapsto \hat{g}_1(t)\hat{g}_2(s)$ is the smooth singular 2-chain with vertices 
in $\mathbf{1}$, $g_1$ and $g_1g_2$ and bounded by the paths $\hat{g}_1$, $g_1\hat{g}_2$ and $\hat{g}_1\hat{g}_2$. Here 
$\Delta_2 = \{(t,s)\in \mathbb{R}^2 \ | \ 0\leq s \leq t \leq 1 \}$ is the 2-simplex in the plane. \begin{figure}[htp]
\centering
\includegraphics[width=1.7in]{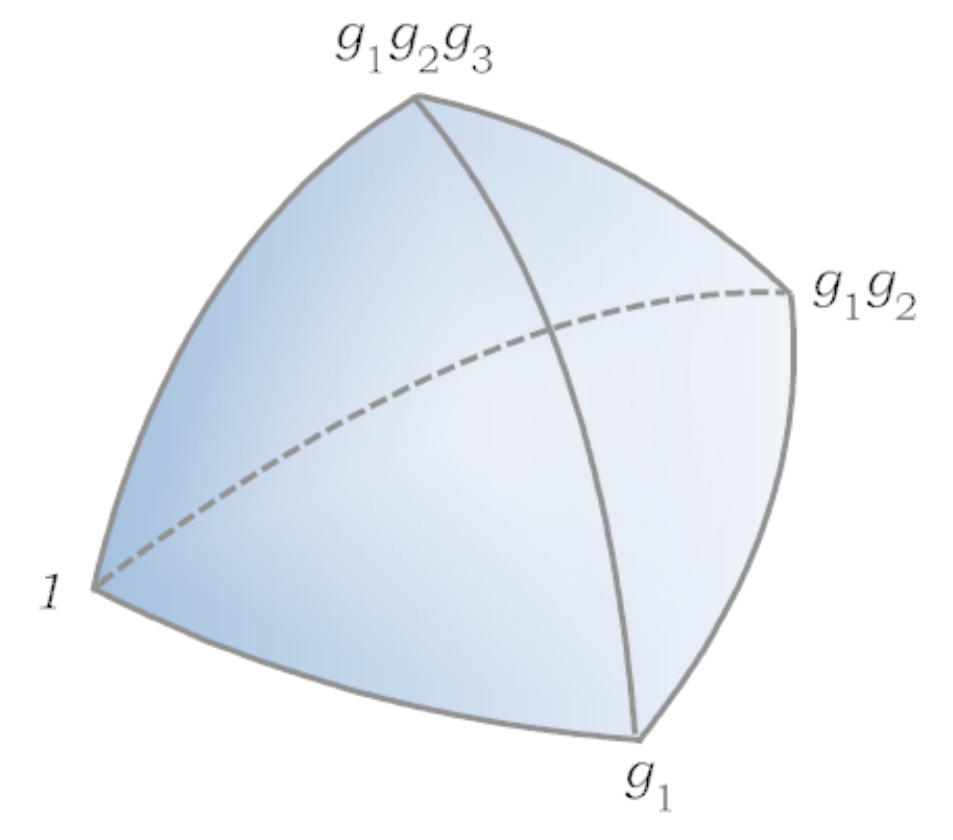}
\caption{Domain of integration in \eqref{2-cocycle}.}\label{Figure1}
\end{figure}

\noindent The 2-cocycle identity 
\begin{equation}\label{2-cocycle}
	(\delta_2\gamma)(\hat{g}_1, \hat{g}_2,\hat{g}_3) = 
\end{equation} 
\[=\hat{g}_1.\gamma(\hat{g}_2,\hat{g}_3) - \gamma(\hat{g}_1\hat{g}_2,\hat{g}_3 )+ \gamma(\hat{g}_1, \hat{g}_2\hat{g}_3) -
	\gamma(\hat{g}_1,\hat{g}_2)= 0 \ \ {\rm mod} \ \Gamma \]

\noindent is satisfied by \eqref{condition} since the regions of integration form a spherical 2-cycle; see Figure 1. 
The face not joining to $\mathbf{1}$ is the left translation by $g_1$ of the domain of integration 
\begin{equation*}
	\hat{g}_1.\gamma(\hat{g}_2,\hat{g}_3) = \int_\sigma g_1.\omega^{eq} = \int_\sigma L_{g_1}^* \omega^{eq} = 
	\int_{L_{g_1} \sigma} \omega^{eq}
\end{equation*}
where we have used the $G$-equivariance of $\omega^{eq}$. Thus, we conclude that the multiplication is associative. 
To see that it is well-defined, i.e. independent of the representatives, a straightforward calculation leads to 

\begin{equation}\label{compatibility 1}
	\int_{\pi[\hat{g}_1\hat{g}_2,\hat{g}'_1\hat{g}_2]}\omega^{eq}- \int_{\pi[\hat{g}_1,\hat{g}'_1]}\omega^{eq} - 
	\gamma(\hat{g}'_1,\hat{g}_2)+\gamma(\hat{g}_1,\hat{g}_2) = 0 \ \ {\rm mod} \ \Gamma \ ,
\end{equation}

\begin{equation}\label{compatibility 2}
	\int_{\pi[\hat{g}_1\hat{g}_2,\hat{g}_1\hat{g}'_2]}\omega^{eq}- \int_{\pi[\hat{g}_2,\hat{g}'_2]}g_1.\omega^{eq} - 
	\gamma(\hat{g}_1,\hat{g}'_2)+\gamma(\hat{g}_1,\hat{g}_2) = 0 \ \ {\rm mod} \ \Gamma \ ,
\end{equation}

\begin{equation}\label{compatibility 3}
	\int_{\pi[\hat{g}_1\hat{g}_2,\hat{g}'_1\hat{g}'_2]}\omega^{eq} - \int_{\pi[\hat{g}_1,\hat{g}'_1]}\omega^{eq} 
	- \int_{\pi[\hat{g}_2,\hat{g}'_2]}g_1.\omega^{eq} -\gamma(\hat{g}'_1,\hat{g}'_2)+\gamma(\hat{g}_1,\hat{g}_2) = 0 
	\ \ {\rm mod} \ \Gamma \ .
\end{equation} 	
Again the regions of integration form closed 2-dimensional surfaces in $G$, depicted in Figure 2. The label on each 
face refers to the corresponding term in the expressions above, numbered from left to right.
\begin{figure}[htp]
\centering
\includegraphics[width=4.8in]{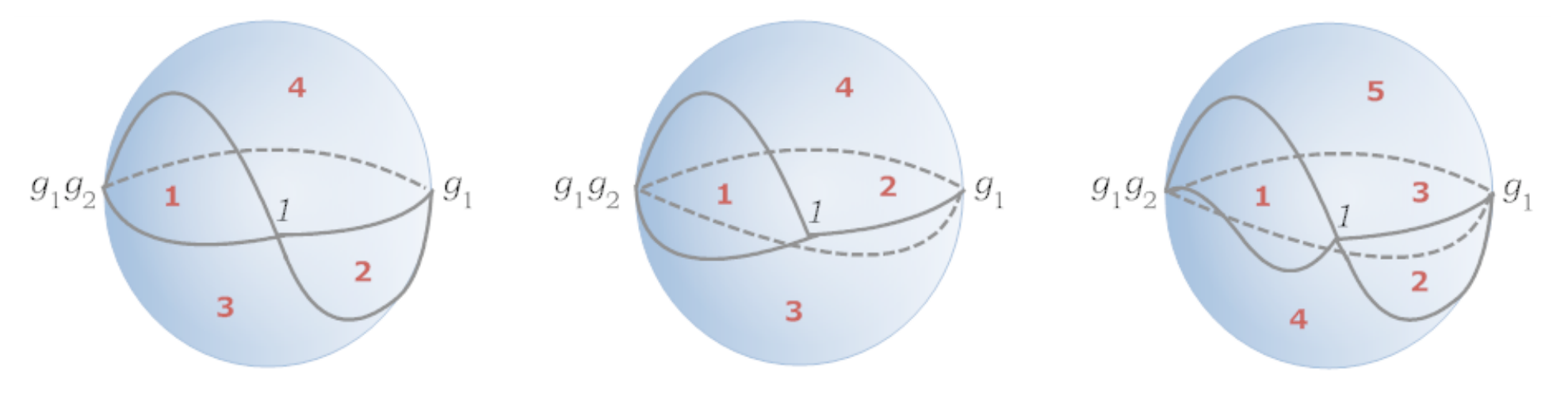}
\caption{Domain of integration in \eqref{compatibility 1}, \eqref{compatibility 2} and \eqref{compatibility 3} respectively.}\label{Figure2}
\end{figure}

Next let us calculate the Lie algebra cocycle. We use the equivariance property to evaluate $\omega^{eq}$ only 
at the identity 

\[	\int_\sigma \omega^{eq} = \int_{\Delta_2}\omega^{eq}(\sigma(t,s))\left(\sigma_*\frac{\partial}{\partial t}, 
	\sigma_*\frac{\partial}{\partial s}\right)dt\wedge ds \]

\[= \int_{\Delta_2}\omega^{eq}(\hat{g}_1(t)\hat{g}_2(s))
	\left(\frac{d\hat{g}_1(t)}{dt}\hat{g}_2(s), \hat{g}_1(t)\frac{d\hat{g}_2(s)}{ds}\right)dtds \] 
	
\[= (g_1 g_2).\int_{\Delta_2}\omega^{eq}(\mathbf{1})\left(\hat{g}^{-1}_2(s)\hat{g}^{-1}_1(t)
	\frac{d\hat{g}_1(t)}{dt}\hat{g}_2(s), \hat{g}^{-1}_2(s)\frac{d\hat{g}_2(s)}{ds}\right)dtds \ .
\]

\noindent For any two elements $X,Y \in \mathcal{P}\mathfrak{g}$, we have smooth curves $c_1(\tau) = \exp(\tau X(t))$ and 
$c_2(\sigma)=\exp(\sigma Y(s))$ in $\mathcal{P}G$, where $\exp$ is defined pointwise by the exponential map of $G$. 
The Lie algebra cocycle is given by

\[(D_2e^\gamma)(X,Y) = \frac{\partial^2}{\partial\sigma\partial\tau}\Big(e^{\gamma(c_1(\tau),c_2(\sigma))}
	e^{-\gamma(c_2(\sigma), c_1(\tau))}\Big)\Big|_{\tau=\sigma=0} \] 
	
\[= \frac{\partial^2}{\partial\sigma\partial\tau}\Big\{\big(\exp(\tau X(1))\exp(\sigma Y(1))\big).\int_{\Delta_2}\omega^{eq}(\mathbf{1})\left(\tau \frac{dX}{dt}, \sigma \frac{dY}{ds}\right)dtds\]

\[ -\big(\exp(\sigma Y(1))\exp(\tau X(1))\big).\int_{\Delta_2}\omega^{eq}(\mathbf{1})\left(\sigma \frac{dY}{ds},\tau \frac{dX}{dt}\right)dtds\Big\}\Big|_{\tau=\sigma=0}\]

\[= \int_{0\leq s \leq t \leq 1} \omega^{eq}(\mathbf{1})\left(\frac{dX}{dt}, \frac{dY}{ds}\right)dtds + 
	\int_{0\leq t \leq s \leq 1} \omega^{eq}(\mathbf{1})\left(\frac{dX}{dt},\frac{dY}{ds}\right)dtds \]

\[= \int_{0\leq s, t \leq 1}\omega^{eq}(\mathbf{1})\left(\frac{dX}{dt}, \frac{dY}{ds}\right)dtds = 
	\omega^{eq}(\mathbf{1})(X(1),Y(1)) = \omega(X,Y)\]

\noindent where we have used the antisymmetry of $\omega^{eq}$ and the fact that $\omega^{eq}(\mathbf{1})(X(0),Y(0)) = 0$. 
Thus, we have verified that $\gamma$ induces a well-defined group multiplication and the correct cocycle at 
the Lie algebra level. The Lie group extension corresponding to $\mathfrak{\hat{g}}=\mathfrak{g}\oplus_\omega
\mathfrak{a}$ is the principal $\widehat{A}$-bundle $ \widehat{G} = \mathcal{P}G\times_\gamma A / \sim \ \to G$ 
with the projection $[(\hat{g},a)] \mapsto \hat{g}(1)$. If $U \subset G$ is a smoothly contractible open identity 
neighborhood, then $\sigma:U\to \widehat{G}, \ g \mapsto [(\hat{g},\mathbf{1}_A)]$ defines a smooth local section, 
where $\hat{g}$ is the path joining $\mathbf{1}$ and $g$ by a smooth contraction. The fiber 
$\widehat{A} = \pi_1(G)\times_\gamma A$ is a central extension of $\pi_1(G)$ by $A$. This can be understood by the 
following argument. If $\pi_1(G)\hookrightarrow \tilde{G}\twoheadrightarrow G$ denotes the universal covering, 
then the same construction for $\tilde{G}$ gives rise to the diagram
\begin{equation}\label{diagram} 
\xymatrix{ 
{A}\ar[r]^-{i_1}\ar[d] 
&{\mathcal{P}\tilde{G}\times_\gamma A / \sim}\ar[r]^-{p_1}\ar[d]^{\phi}&{\tilde{G}}\ar[d]^{q}\\ 
{\widehat A}\ar[r]^-{i_2} &{\mathcal{P}G\times_\gamma A / \sim}\ar[r]^-{p_2}&{G}} 
\end{equation}
Restriction to the subgroup $\pi_1(G) \subset \tilde{G}$ induces a central extension of $\pi_1(G)$,
\begin{equation*}
	A \xrightarrow{i_1} \pi_1(G)\times_\gamma A \xrightarrow{p_1} \pi_1(G) \ ,
\end{equation*}
since $\pi_1(G)$ is discrete and acts trivially on $A$. Finally we have $\widehat{A} = {\rm ker \ } p_2 = {\rm ker \ } 
q \circ p_1 = \pi_1(G)\times_\gamma A$. The right action of the structure group $\widehat{A}$ is given by 
$[(\hat{g},a)].[(\eta, a')] = [(\hat{g}\eta, a(\hat{g}.a')e^{\gamma(\hat{g}, \eta)})]$. 

To see that condition \eqref{condition} is not only sufficient but necessary, let $A\hookrightarrow\widehat{G}\twoheadrightarrow G$ 
be a Lie group extension corresponding to $\mathfrak{\hat{g}} = \mathfrak{g} \oplus_\omega \mathfrak{a}$. There
is a connection form on this principal bundle whose curvature is the \textit{right} $G$-equivariant 2-form $\omega^{eq}_R \in \Omega^2(G,\mathfrak{a})$, 
\begin{equation*}
	\omega^{eq}_R(g)(R_{g*}X, R_{g*}Y) = (R_{g}^*\omega^{eq})(\mathbf{1})(X,Y) = g^{-1}.\omega(X,Y) \ \ \ \ \ \forall X,Y \in 
	\mathfrak{g} \ ,
\end{equation*}
where $R_g:G\to G, h\mapsto hg$ denotes the right translation map. Indeed, if 
$pr_\mathfrak{a}: \mathfrak{g}\oplus_\omega \mathfrak{a} \to \mathfrak{a}$ denotes the projection onto the ideal $\mathfrak{a}$, 
then there is a canonical connection 1-form on $\widehat{G}$ given by
\begin{equation}\label{connection}
	\alpha = \hat g^{-1} pr_\mathfrak{a}\big(d\hat g \hat g^{-1}\big) \hat g \in \Omega^1(\widehat{G},\mathfrak{a})
\end{equation}
with the curvature $\Omega = d\alpha + \frac{1}{2}[\alpha,\alpha] = d\alpha$, where the commutator term vanishes since 
$\mathfrak{a}$ is abelian. By the Maurer-Cartan equation it follows that  
\begin{eqnarray*}
	d\alpha &=&  d\hat g^{-1} pr_\mathfrak{a}\big(d\hat g \hat g^{-1}\big) \hat g - \hat g^{-1} pr_\mathfrak{a}\big(d\hat g d\hat g^{-1}\big) \hat g - \hat g^{-1} pr_\mathfrak{a}\big(d\hat g \hat g^{-1}\big) d\hat g \\
&=& -[\hat g^{-1} d\hat g, \hat g^{-1} pr_\mathfrak{a}\big(d\hat g \hat g^{-1}\big) \hat g] +  \frac 12 \hat g^{-1} pr_\mathfrak{a} \big([d\hat g \hat g^{-1}, d\hat g \hat g^{-1}]\big) \hat g \ .
\end{eqnarray*}
Since $\mathfrak{\hat{g}} = \mathfrak{g} \oplus_\omega \mathfrak{a}$ is a direct sum, we may write $d\hat g \hat g^{-1} = \theta_\mathfrak{g}+\theta_\mathfrak{a}$ as a sum of two 1-forms, where $\theta_\mathfrak{g}$ and $\theta_\mathfrak{a}$ only have components along the $\mathfrak{g}$ and $\mathfrak{a}$ directions respectively. We have 	
\begin{eqnarray*}d\alpha &=& -[\hat g^{-1} \theta_\mathfrak{g} \hat g, \hat g^{-1} \theta_\mathfrak{a} \hat g] +   \hat g^{-1} pr_\mathfrak{a}\big(\theta_\mathfrak{g}^2 +\theta_\mathfrak{g}\theta_\mathfrak{a}+\theta_\mathfrak{a}\theta_\mathfrak{g}+\theta_\mathfrak{a}^2\big) \hat g \\
 &=& - \hat g^{-1} pr_\mathfrak{a}\big(\theta_\mathfrak{g}\theta_\mathfrak{a}+\theta_\mathfrak{a}\theta_\mathfrak{g}\big) \hat g + \hat g^{-1} \omega(dg g^{-1}, dg g^{-1}) \hat g + \hat g^{-1} pr_\mathfrak{a}\big(\theta_\mathfrak{g}\theta_\mathfrak{a}+\theta_\mathfrak{a}\theta_\mathfrak{g}\big) \hat g \\
&=& g^{-1}.\omega(dg g^{-1}, dg g^{-1}) = \omega^{eq}_R(g) \ .
\end{eqnarray*}
Note that in case of central extensions, there exists a different connection form $\tilde \alpha = pr_\mathfrak{a}\big(\hat g^{-1} d\hat g \big)$ on $\widehat G$ whose curvature is precisely the left invariant 2-form determined by $\omega$. However, when $G$ acts non-trivially on $A$, the 1-form $\tilde \alpha$ is no longer a connection since it is not invariant under right translation by elements in $A$ and we must therefore use (\ref{connection}).

Equipped with a connection, we can define the horizontal lift of a curve on the base space and subsequently the 
notion of parallel transport. In particular, the holonomy around a smooth contractible loop $\eta:\mathbb T \to G$ is 
given by
\begin{equation*}
	{\rm hol(\alpha, \eta)} = e^{ \int_\eta \alpha} = e^{ \int_\pi \omega_R^{eq} }
\end{equation*}
where $\mathbb T =  \mathbb R / \mathbb Z$ denotes the unit circle and $\pi$ is a surface enclosed by $\eta$. 
The arbitrariness in the choice of this surface leads as before to the requirement $\int_{[c]} \omega_R^{eq} \in \Gamma$ for all spherical cycles $[c] \in H_2(G)$.
The integrality condition $\int_{[c]} \omega^{eq} \in \Gamma$ then follows by the diffeomorphism $g\mapsto g^{-1}$,
which relates the right and left $G$-equivariant forms.
Moreover, since spherical cycles factor through to spheres, the condition can be formulated 
equivalently in terms of $\pi_2(G) \cong H_2(\tilde G)$. Thus, we are led to the result:

\begin{theorem}[Integrability criterion] 
	Let $G$ be a connected Lie group and $A$ a smooth $G$-module of the form $\mathfrak{a}/\Gamma$ for some 
	discrete subgroup $\Gamma \subseteq \mathfrak{a}$. The abelian Lie algebra extension 
	$\hat{\mathfrak{g}} = \mathfrak{g}\oplus_\omega \mathfrak{a}$ integrates to a Lie group extension 
	$\mathbf{1} \to \widehat{A}\to \widehat{G}\to G \to \mathbf{1}$ if and only if the image of the period
	homomorphism
	\begin{equation*}
		per_\omega:\pi_2(G) \to \mathfrak{a}, \ \ \  [\sigma]  \mapsto \int_{[\sigma]} \omega^{eq} 
	\end{equation*}
	is contained in $\Gamma$, where $\widehat{A} = \pi_1(G)\times_\gamma A$ is a central extension of 
	$\pi_1(G)$ by $A$. 
\end{theorem}
When $G$ is simply connected, we have $\widehat{A} = A$ and the condition coincides with that found in \cite{N}. 
More generally, whenever the central extension $\widehat{A}$ splits, we may mod out by $\pi_1(G)$ and the Lie 
algebra extension lifts to $A \hookrightarrow \widehat{G} \twoheadrightarrow G$.

\begin{example}
	Consider the two-dimensional torus $G = \mathbb T^2 = \mathbb R^2 / \mathbb Z^2$ and the left invariant 2-form 
	corresponding to the symplectic form 
	\begin{equation*} \omega = dx\wedge dy\end{equation*}
	on the universal covering group $\tilde{G} = \mathbb R^2$. 
	This form is invariant under the $\pi_1(G) = \mathbb Z^2$ action and since $\pi_2(G)=0$, all its periods are integral.
	The central extension $\widehat{\mathfrak{g}}=\mathbb R^2\oplus_\omega \mathbb R$ of the Lie algebra $\mathfrak{g}=\mathbb R^2$ 
	by the trivial $\mathfrak{g}$-module $\mathbb R$ defines the Heisenberg group, and according to Theorem 4.1 it lifts to a central 
	extension 
	\begin{equation*}
		\mathbb{Z}^2\times_\gamma \mathbb T \hookrightarrow \widehat{\mathbb T^2}\twoheadrightarrow \mathbb T^2 \ .
	\end{equation*}	
	In this case, the central extension
	\begin{equation*}
		\mathbb T \hookrightarrow \mathbb{Z}^2\times_\gamma \mathbb T \twoheadrightarrow \mathbb{Z}^2
	\end{equation*} 
	does not split and can be characterized as follows. The fundamental
	group of the torus is the integral lattice in $\mathbb R^2$, and to 
	each point in the lattice one can associate the path connecting it to 
	the origin, for example, the straight line. For a pair of straight lines 
	$\hat{g}_1(t)$ and $\hat{g}_2(t)$, the 2-cocycle $\gamma$ is computed as the 
	integral of the 2-form $\omega$ over a two-dimensional surface determined by the paths. 
	The end result is a 2-cocycle $e^{i \beta(\textbf{p},\textbf{q})}$ for $(\textbf{p},\textbf{q})\in \mathbb Z^2\times \mathbb Z^2$, where $\beta$ is a real antisymmetric 
	bilinear form on the integral lattice. The value of the form $\beta$ is simply
	the area of the triangle with vertices at $\textbf{0}$, $\textbf{p}$ and $\textbf{p}+\textbf{q}$ (up to a normalization 
	of the symplectic form).
	
	Furthermore, in the universal covering group $\mathbb R^2$, any two paths with the same endpoints are homotopic, and 
	hence there is no discrete modification of the fibre $\mathbb T$ in the corresponding central extension,
	\begin{equation*}
		\mathbb T \hookrightarrow \widehat{\mathbb R^2}\twoheadrightarrow \mathbb R^2 \ .
	\end{equation*}
\end{example}

\vspace{30pt}
\begin{remark} Subsequent to the submission of this paper to the preprint archive \texttt{arXiv} (math/0611431), a similar approach to abelian Lie group extensions using path groups appeared in (math/0703342) \cite{V}. Vizman's idea is to embed the group of contractible loops $\Omega_0G \subset \Omega G$ as a normal subgroup in $\mathcal{P}G\times_\gamma A$ via the graph of the smooth map
\begin{equation*}
		\lambda:\Omega_0G\to A, \ \ g\mapsto \left(e^{\int_{\hat g}\omega^{eq}}\right)^{-1} \ ,
\end{equation*}
where $\hat g$ is an arbitrary path in $\Omega_0G$, interpreted as a map $\hat g: [0,1]^2\to G$. The quotient group $\mathcal{P}G\times_\gamma A / \text{Graph}(\lambda)$ defines an abelian extension of the universal covering group $\tilde{G}$ by $A$. This method is closer in spirit to the constructions in \cite{LMNS} and \cite{M3}. Vizman has cited this paper and explained in Remark 8 \cite{V} that the extension $\mathcal{P}G\times_\gamma A / \text{Graph}(\lambda)$ is canonically isomorphic to the Lie group $\mathcal{P}\tilde G\times_\gamma A /\sim$ in diagram \eqref{diagram}. 
\end{remark} 

\section*{Acknowledgements}
\noindent The author thanks Jouko Mickelsson and Karl-Hermann Neeb for helpful remarks and
is grateful to the referees for their thorough reviews and suggestions, which greatly improved the presentation of this paper.
\bibliographystyle{amsplain}

\begin{thebibliography}{10}

\bibitem {CM} A. L. Carey and M. K. Murray, \textit{String Structures and the Path Fibration of a Group}, Comm. Math. Phys. \textbf{141}, no. 3 (1991), pp. 441--452.

\bibitem {E} S. Eilenberg, \textit{On Spherical Cycles}, Bull. Amer. Math. Soc. \textbf{47}, no. 6 (1941), pp. 432--434.

\bibitem {G} H. Gl\"ockner, \textit{Fundamental Problems in the Theory of Infinite-Dimensional Lie Groups}, \\ J. Geom. Symm. Phys. \textbf{5} (2006), pp. 24--35.

\bibitem {K1} A. Kriegl and P. W. Michor, \textit{Regular Infinite Dimensional Lie Groups}, J. Lie Theory \textbf{7} (1997), pp. 61--99.

\bibitem {K2} A. Kriegl and P. W. Michor, \textit{The Convenient Setting of Global Analysis}, Amer. Math. Soc., Providence, (1997).

\bibitem {LMNS}{A. Losev, G. Moore, N. Nekrasov and S. Shatashvili}, \textit{Central Extensions of Gauge Groups 
Revisited}, Sel. Math., New Ser. 4(1) (1998), pp. 117--123.

\bibitem {M1} P. W. Michor and  J. Teichmann, \textit{Description of Infinite-Dimensional Abelian Regular Lie Groups}, J. Lie Theory \textbf{9} (1999), pp. 487--489.

\bibitem {M2} J. Mickelsson, \textit{Current Algebras and Groups}, Plenum Press, New York (1989).

\bibitem {M3} J. Mickelsson, \textit{Kac-Moody Groups, Topology of the Dirac Determinant Bundle and Fermionization}, Comm. Math. Phys. \textbf{110}, no. 2 (1987), pp. 173--183.

\bibitem {M4} M. K. Murray, \textit{Another Construction of the Central Extension of the Loop Group}, Comm. Math. Phys. \textbf{116}, no. 1 (1988), pp. 73--80.

\bibitem {M5} J. Milnor, \textit{Remarks on Infinite-Dimensional Lie Groups}, ``Relativit\'e, Groupes et Topologie II'',
B. DeWitt and R. Stora (Eds.), North-Holland, Amsterdam (1983), pp. 1007--1057.

\bibitem {N} K. H. Neeb, \textit{Abelian Extensions of Infinite-Dimensional Lie Groups}, Mathematical Works. Part XV. Luxembourg: Universit\'{e} du Luxembourg, S\'{e}minaire de Math\'{e}matique (2004), \\ pp. 69--194.

\bibitem {PS} A. Pressley and G. Segal, \textit{Loop Groups}, Clarendon Press, Oxford (1986).

\bibitem {V} C. Vizman, \textit{The Path Group Construction of Lie Group Extensions}, J. Geom. Phys. \textbf{58}, \\ no. 7 (2008), pp. 860--873.

\end{thebibliography}

\end{document}